\documentclass{amsart}

\usepackage{amssymb}
\usepackage{amsmath}

\usepackage[english,francais]{babel}

 \newtheorem{theorem}{Theorem}[section]

 \theoremstyle{definition}
 \newtheorem{definition}[theorem]{Definition}
 \theoremstyle{remark}

 \numberwithin{equation}{section}

\DeclareMathOperator{\tr}{\mathrm{trace}}

\def\R {\mathbb{R}}       
\def\N {\mathbb{N}}       

\def\Rn {\mathbb{R}^{n}}

\def\mr{\mathrm}

\def\cL{{\mathcal L}}

\def\cR{{\mathcal R}}
\def\cS{{\mathcal S}}

\def\1{\lambda}
\def\2{\Sigma_{2}}

\def\<{{\langle}}
\def\>{{\rangle}}

\def\H{{\mathcal H}}

\def\Rn{{\R^n}}

\def\Gh{\widehat{G}}
\begin{document}

\title[Lower bounds for operators on graded Lie groups]
{Lower bounds for operators on graded Lie groups}
\author[V. Fischer]{V\'eronique Fischer}

\address{Universita degli studi di Padova, DMMMSA, Via Trieste 63,
           35121 Padova, Italy}

\email{fischer@dmsa.unipd.it}

\author[M. Ruzhansky]{Michael Ruzhansky}
\address{180 Queen's Gate, Department of Mathematics, Imperial College London, London, 
SW7 2AZ, United Kingdom}
\email{m.ruzhansky@imperial.ac.uk}

\subjclass{Primary 35S05; Secondary 43A80}

\keywords{harmonic analysis,
nilpotent Lie groups, 
pseudo-differential operators}

\date{February, 2013}

\maketitle
\begin{abstract}
\selectlanguage{english}
In this note we present a symbolic pseudo-differential calculus 
on any graded (nilpotent) Lie group
and, as an application, a version of the sharp G{\aa}rding inequality. 
As a corollary, we obtain lower bounds for
positive Rockland operators with variable coefficients as well as 
their Schwartz-hypoellipticity.
\end{abstract}

\vskip 0.5\baselineskip

\section{Introduction}
\label{Intro}

In this note we present a symbolic pseudo-differential calculus on any
graded Lie group.
As applications, we obtain a version of the sharp G{\aa}rding inequality and
results on the Schwartz hypoellipticity for operators in this context.

In the usual Euclidean setting,
the positivity of the full symbol is required
 for the sharp G{\aa}rding inequality
as well as, for instance, for the Fefferman-Phong inequality.
This contrasts with
 the G{\aa}rding  or Melin-H\"ormander inequalities, for example,
where knowing the principal (or subprincipal) symbol is sufficient.
Thus the latter inequalities can be proved on manifolds using the standard
H\"ormander theory of pseudo-differential operators together with 
the usual Kohn-Nirenberg quantisation on $\Rn$.
Since the geometric control of the full
symbol of an operator  
is impossible with these tools in general,
the study of sharp G{\aa}ring inequalities appears limited.
However, 
the sharp G{\aa}rding inequality on compact Lie groups 
was recently established in \cite{RTg}.
This approach uses the notion of a full matrix-valued
symbol defined in terms of the representation theory of the group.
In this note we explain how we follow the same strategy in the case of the Heisenberg group or more general nilpotent
Lie groups. 

The pseudo-differential calculus that we define is different from 
the several  pseudo-differential calculi already developed 
on the Heisenberg group: see e.g. Taylor
\cite{Tnma} and \cite{Tnha} for symbol classes coming from standard
H\"ormander classes though the exponential mapping;
Bahouri, Fermanian-Kammerer and Gallagher \cite{BFG} 
for classes defined in terms of explicit formulae
coming from the Schr\"odinger representation; or
Beals and Greiner \cite{BG} or Ponge \cite{P} for different types of analysis on
the Heisenberg manifolds. On more general nilpotent groups,
Christ, Geller, Glowacki and Polin \cite{CGGP} proposed an approach to 
pseudo-differential operators, however based on the properties of kernels and not on a symbolic calculus.
Following Ruzhansky and Turunen \cite{RTi} and \cite{RTb}, 
we define symbol classes directly
on the group. As such, our approach can be extended for general
graded nilpotent Lie groups, and by developing the symbolic
calculus and the Friedrichs approximation on the group,
we obtain the corresponding sharp G{\aa}rding inequality
and Schwartz hypoellipticity results. 

While the symbol classes in \cite{RTb} are based on the spectral
theory of the Laplace-Beltrami operator, here, it is not available and it
becomes natural to use the
sub-Laplacian on stratified groups or more general Rockland
operators on graded groups. 
Moreover, surpassing \cite{RTg},
since a dilation structure is present, we establish the
sharp G{\aa}rding inequality for suitable $(\rho,\delta)$ classes
of operators by establishing a Calder\'on--Vaillancourt type
theorem in this context. This is, in fact, the best known lower bound available
in the $(\rho,\delta)$-setting already on $\Rn$.

The appearing operators are Calder\'on--Zygmund in the sense of
Coifman and Weiss \cite[ch. III]{CW}, so that $L^{p}$ results follow as well.
In Section \ref{prel} we fix the notation concerning Lie groups that
we are working on. In Section \ref{res} we formulate the results.

\section{Preliminaries}
\label{prel}

Let us first briefly recall the necessary notions and set some notation.
In general, we will be concerned with graded Lie groups $G$
which means that $G$ is a connected and simply connected 
Lie group (of step $s>1$) with the gradation of its Lie algebra $\mathfrak g$ given by
$\mathfrak g=
\underset {\ell=1}  {\overset \infty \oplus} \mathfrak g_{\ell}$
with
$[\mathfrak g_{\ell},\mathfrak g_{\ell'}]\subset\mathfrak g_{\ell+\ell'}$
for any $\ell,\ell'\in \N$, 
where the $\mathfrak g_{\ell}$, $\ell=1,2,\ldots$ 
are vector subspaces of $\mathfrak g$,
almost all equal to $\{0\}$.
This implies that the group $G$ is nilpotent.
If the whole $\mathfrak g$ is generated
by $\mathfrak g_{1}$ in this way, the group $G$ is said to be stratified.

Let $\{X_1,\ldots X_{n_1}\}$ be a basis of $\mathfrak g_1$ (this basis is possibly reduced to $\{0\}$), let
$\{X_{n_1+1},\ldots,  X_{n_1+n_2}\}$ a basis of $\mathfrak g_2$
and so on, so that we obtain a basis 
$X_1,\ldots, X_n$ of $\mathfrak g$ adapted to the gradation.
Via the exponential mapping $\exp_G : \mathfrak g \to G$, we   identify 
the points $(x_{1},\ldots,x_n)\in \R^n$ 
 with the points  $x=\exp_G(x_{1}X_1+\cdots+x_n X_n)$ in $G$.
This leads to a corresponding Lebesgue measure on $\mathfrak g$ and the Haar measure $dx$ on the group $G$.
We define
 the  spaces of Schwartz functions  $\cS(G)$ and tempered distributions $\cS'(G)$ of the group $G$  as  those on $\R^n$.
The coordinate function $x=(x_1,\ldots,x_n)\in G\mapsto x_j \in \R$
is denoted by $x_j$.
More generally we define for every multi-index $\alpha\in \N_0^n$,
$x^\alpha:=x_1^{\alpha_1} x_2 ^{\alpha_2}\ldots x_{n}^{\alpha_n}$, 
as a function on $G$.
Similarly we set
$X^{\alpha}=X_1^{\alpha_1}X_2^{\alpha_2}\cdots
X_{n}^{\alpha_n}$ in the universal enveloping Lie algebra of $\mathfrak g$.

For any $r>0$, 
we define the  linear mapping $D_r:\mathfrak g\to \mathfrak g$ by
$D_r X=r^\ell X$ for every $X\in \mathfrak g_\ell$, $\ell\in \N$.
Then  the Lie algebra $\mathfrak g$ is endowed 
with the family of dilations  $\{D_r, r>0\}$
and becomes a homogeneous Lie algebra in the sense of \cite{FS}.
 The weights of the dilations are the integers $\upsilon_1,\ldots, \upsilon_n$ given by $D_r X_j =r^{\upsilon_j} X_j$, $j=1,\ldots, n$.
 The associated group dilations are defined by
$$
r\cdot x
:=(r^{\upsilon_1} x_{1},r^{\upsilon_2}x_{2},\ldots,r^{\upsilon_n}x_{n}),
\quad x=(x_{1},\ldots,x_n)\in G, \ r>0.
$$
In a canonical way  this leads to the notions of homogeneity for functions and operators.
For instance
the degree of homogeneity of $x^\alpha$ and $X^\alpha$,
viewed respectively as a function and a differential operator on $G$, is 
$[\alpha]=\sum_j \upsilon_j\alpha_{j}$.
Indeed, let us recall 
that a vector of $\mathfrak g$ defines a left-invariant vector field on $G$ 
and more generally 
that the universal enveloping Lie algebra of $\mathfrak g$ 
is isomorphic with the left-invariant differential operators; 
we keep the same notation for the vectors and the corresponding operators. 

The dimension of $G$ is
$n=\sum_\ell n_\ell$ while its homogeneous dimension is
$Q=\sum_\ell  \ell n_\ell = \upsilon_1+\upsilon_2+\ldots+\upsilon_n$.

We denote by $\Gh$ the set of equivalence classes of (continuous) irreducible
unitary representations of $G$.
We will often identify a representation of $G$ with its equivalence class.
We will also keep the same notation for the corresponding infinitesimal representation. 
For $\pi\in \Gh$, 
 we denote by $\H_{\pi}$ the representation space of $\pi$ and by
$\H_{\pi}^{\infty}$ its subspace of smooth vectors.
For $f\in L^{1}(G)$, we define its Fourier transform at
$\pi\in\Gh$ by
$\widehat{f}(\pi)=\int_{G} f(g) \pi(g)^{*} dg$, 
with the integral understood in the
Bochner sense. 
Denoting by $\mu$ the Plancherel measure on $\Gh$,
the inverse Fourier formula holds:
$$
f(g) = \int_{\Gh} \tr \left(\pi(g) \widehat{f}(\pi) \right)
d\mu(\pi)
\quad\mbox{when}\quad
\int_{\Gh} \tr \left|\widehat{f}(\pi) \right|
d\mu(\pi)<\infty
\ .
$$

Let $\cR$ be a positive (left) Rockland operator on $G$; 
this means that $\cR$
is a left-invariant differential operator, homogeneous of degree $\nu$ necessarily even, positive in
the operator sense, and such
that for every non-trivial $\pi\in\Gh$ the operator
$\pi(\cR)$ is injective on $\H_\pi^\infty$. 
The operator
$\cR$ admits an essentially self-adjoint
extension on $C_{0}^{\infty}(G)$ (see \cite{FS}), and we will still denote this extension by $\cR$.
Examples of such operators are given 
in the stratified case by $\cR=-\cL$ 
where $\cL=\sum_{1\leq j\leq n_{1}} X_i^{2}$
is a Kohn-sub-Laplacian,
and in the graded case by the operators
$$
\sum_{1\leq j \leq n_j }
(-1)^{\frac{\nu_o}{\upsilon_j}} X_j^{2\frac{\nu_o} {\upsilon_j}} 
\quad\mbox{and}\quad
\sum_{1\leq j \leq n_j }
X_j^{4\frac{\nu_o} {\upsilon_j}} 
, $$
where $\nu_o$ denotes some common multiple of $\upsilon_1,\ldots,\upsilon_n$.
In fact our class of operators do not depend on the choice of such an operator $\cR$.

\section{Results}
\label{res}

We aim at defining the symbol classes in terms of the operators $\cR$ as above. 

\begin{definition}
A \emph{symbol} is a family of operators
$\sigma=\{\sigma(x,\pi):\; x\in G, \; [\pi]\in\Gh\}$,
such that
\begin{itemize}
\item[1.] for each $x\in G$, the family $\{\sigma(x,\pi), \pi\in \hat G\}$ is a $\mu$-measurable field of
operators $\H_{\pi}^{\infty}\to\H_{\pi}$;
\item[2.] there exist two constants $\gamma_{1},\gamma_{2}\in\R$ such that for every $x\in G$,
the operator 
$\pi(I+\cR)^{\gamma_{1}}\sigma(x,\pi)\pi(I+\cR)^{\gamma_{2}}$ is
bounded on $\H_\pi$ uniformly in $\pi\in \hat G$;
\item[3.] for any $\pi\in\Gh$ and any $u,v\in\H_{\pi}$, the scalar function
$x\mapsto (\sigma(x,\pi)u,v)_{\H_{\pi}}$ is smooth on $G$.
\end{itemize}
\end{definition}
Here, the powers $\pi(I+\cR)^{\gamma}$ are defined by the spectral theorem 
for the positive operator $\pi(\cR)$.
The existence of $\gamma_{1},\gamma_{2}$ in the second condition is used to guarantee that the following formula makes sense:
$$
Tf(x)=\int_{\Gh} \tr\left(\pi(x)\sigma(x,\pi)\widehat{f}(\pi)\right) d\mu(\pi)
\quad, \quad f\in \cS(G), \, x\in G 
\ ;
$$
indeed such operator $T={\mr{Op}}(\sigma)$ is 
well-defined and continuous $\cS(G)\to\cS'(G)$.
We note that if the operator $T$ is left-invariant, then 
its symbol is independent of $x$.

\begin{definition}
The \emph{difference operators} $\Delta^\alpha$, $\alpha\in \N_0^n$,
are densely defined on the $C^*$-algebra of the group 
via
$$
(\Delta^\alpha \widehat{f})(\pi):=\widehat{(x^\alpha f)}(\pi)
\quad,\quad 
f\in \cS(G)\ .
$$
\end{definition}

Let now $m\in\R$ and $0\leq\delta\leq\rho\leq 1$ with $\delta\not=1$.
\begin{definition}
The symbol class $S^{m}_{\rho,\delta}$ is defined as the set of  symbols
$\sigma$ satisfying for all $\alpha,\beta\in \N_0^n$ and every $\gamma\in\R$:
\begin{equation}\label{EQ:symbols}
\sup_{x\in G, \pi\in\Gh} \| \pi(I+\cR)^{\frac{\rho[\alpha]-m-\delta[\beta]+\gamma}{\nu}}
X_{x}^{\beta}\Delta^{\alpha}\sigma(x,\pi)
\pi(I+\cR)^{-\frac{\gamma}{\nu}}
\|_{op}<\infty.
\end{equation}
(The supremum over $\pi$ is in fact the essential supremum over the Plancherel measure $\mu$.)
\end{definition}
It is easy to see that if $m_{1}\leq m_{2}$,
$\delta_{1}\leq\delta_{2}$ and $\rho_{1}\geq\rho_{2}$, then
$S^{m_{1}}_{\rho_{1},\delta_{1}}\subset S^{m_{2}}_{\rho_{2},\delta_{2}}$.
Furthermore, the expressions in \eqref{EQ:symbols} 
define a Fr\'echet topology on
the linear space $S^{m}_{\rho,\delta}.$ 

In the abelian case,
that is, $\R^n$ endowed with the addition law
and $\cR=-\cL$, $\cL$ being the Laplace operator, 
$S^m_{\rho,\delta}$ boils down easily to the usual H\"ormander class.
However our initial motivation did not come from the abelian case: 
we wanted to define 
the difference operators and the symbol classes 
in analogy with the ones defined in \cite{RTi} on compact Lie groups.
In this case, a definition similar to \eqref{EQ:symbols}
 would formally give the same classes of symbols defined in \cite{RTi} 
since, $\cR=-\cL$,  $\cL$ being the Laplace-Beltrami operator, 
the operator $\pi(I+\cR)$ is scalar.
In our case, i.e. $G$ being a graded non-abelian Lie group, 
the operator $\cR$ is not even central
and the introduction of $\gamma$ in \eqref{EQ:symbols} assures
that $\bigcup_{m\in \R} S^m_{\rho,\delta}$ is an algebra. 

We have the following properties
for the operators classes
$\Psi^{m}_{\rho,\delta}:={\mr{Op}}(S^{m}_{\rho,\delta})$ defined using 
the quantisation procedure $\sigma \mapsto {\mr{Op}}(\sigma)$ described above.

\medskip

\begin{theorem}\label{thm:1}
Let $0\leq \delta\leq \rho\leq 1$.
We have the following properties:
\begin{itemize}
\item[(1)] 
The symbol classes is an algebra of operators $\bigcup_{m\in \R} S^m_{\rho,\delta}$  
stable by taking the adjoint. 
Each vector space $S^m_{\rho,\delta}$ does not depend on the choice of the positive Rockland operator $\cR$.
\item[(2)]
For $\rho\not=0$, the operator class
 $\bigcup_{m\in \R} \Psi^m_{\rho,\delta}$  is an algebra
stable by taking the adjoint.
\item[(3)] 
For any $\alpha\in\N_0^n$, we have $X^{\alpha}\in\Psi^{[\alpha]}_{1,0}.$
\item[(4)] For any positive Rockland operator of homogeneous degree $\nu$, 
we have 
$(I+\cR)^{\frac m \nu}\in \Psi^m_{1,0}$.
\item[(5)]
If $\rho\in [0,1)$
then the operators in $\Psi^0_{\rho,\rho}$ are continuous on $L^2(G)$.
\item[(6)] Let $\rho\not=0$.
The integral kernel $K(x,y)$ of an operator $T\in\Psi^m_{\rho,\delta}$ is
smooth on $(G\times G)\backslash\{(x,y):x=y\}$.
It is of Calder\'on-Zygmund type in the sense of  Coifman and Weiss \cite[ch.III]{CW}. 
It decreases rapidly  as $|xy^{-1}|\to \infty$
(here we have fixed a homogeneous norm $|\cdot|$ on $G$, 
 i.e. a continuous function, homogeneous of degree one and
vanishing only at $0$):
i.e. for any $M>0$ there exists $C_{M}>0$ such that 
$$
|xy^{-1}|\geq 1 \ \Longrightarrow \
|K(x,y)|\leq C_{M}|xy^{-1}|^{-M}
\ .$$ 
 At the diagonal it satisfies
$$
|xy^{-1}|\leq 1 \ \Longrightarrow \
|K(x,y)|\leq C|xy^{-1}|^{-\frac{Q+m}{\rho}}\ .
$$ 
\end{itemize}
\end{theorem}
\medskip

By \cite[ch.III th\'eor\`eme 2.4]{CW} 
Property (6) implies that
the operators in $\Psi^0_{\rho,\delta}$,
  $1\geq \rho\geq \delta\geq 0$, $\rho\not=0$, $\delta\not=1$,
 are continuous on $L^{p}(G)$, $1<p<\infty$.

By Properties (2) and (4), any operator in $\Psi^m_{\rho,\delta}$, 
$1\geq \rho\geq \delta\geq 0$, $\rho\not=0$, $\delta\not=1$,
is continuous on the natural Sobolev spaces (denoted by $L^2_a(G)$) associated with the dilations 
and the loss of derivatives is controlled by the order $m$. 
The Sobolev space $L^{2}_{a}(G)$ is defined as the set 
of tempered distribution $f\in \cS'(G)$ such that 
$(I+\cR)^{\frac{a}{\nu}}f\in L^{2}(G)$ but does not depend on the choice of $\cR$.
These Sobolev spaces enjoy properties similar to the stratified case proved by Folland \cite{F}, in particular for interpolation (see \cite{RF}).

We now give the sharp G{\aa}rding inequality. 
\medskip

\begin{theorem}\label{thm:2}
Let $0\leq\delta\leq\rho\leq 1$, $\rho\not=0$, $\delta\not=1$,
and let $T\in\Psi^{m}_{\rho,\delta}$ with symbol $\sigma=\{\sigma(x,\pi)\}$.
Assume that each $\{\sigma(x,\pi)\}$ of $T$ is non-negative on $\H_\pi$ (in the operator sense).
Assume also that there exists a Rockland operator $\cR$ such that 
each $\sigma(x,\pi)$ commutes with the spectral measure of $\pi(\cR)$ for every
$x\in G$ and almost every $\pi\in\Gh$. 
 Then there exists $C>0$ such that for
every $f\in\cS(G)$ we have
$$
{\mr{Re}} (Tf,f)_{L^{2}(G)}\geq -C\|f\|_{L^{2}_{\frac{m-(\rho-\delta)}{2}}(G)}.
$$
\end{theorem}

The class includes 
\begin{itemize}
\item 
the variable
coefficient Kohn-sub-Laplacians and Rockland operators of the form
$a(x)\cR$, with $a(x)\geq 0$ satisfying
$X^{\alpha}a\in L^{\infty}(G)$ for all $\alpha\in \N_0^n$,
\item the multipliers $\phi(\cR)$ for 
a smooth function $\phi:[0,\infty) \mapsto [0,\infty)$ satisfying
\begin{equation}
\label{eq_cond_mutliplier}
\forall a\in \N_0\quad \exists C=C_a>0\qquad
\forall \lambda\geq 0\quad
|\partial_\lambda^\alpha \phi(\lambda)|\leq C (1+\lambda)^{\frac m \nu -a},
\end{equation}
\item more generally the operators with symbols 
given by $\sigma(x,\pi)=\phi_x(\pi(\cR))$
with $(x,\lambda)\mapsto \phi_x(\lambda)$ being non-negative and smooth on 
$G\times [0,\infty)$,
and $\phi_x$ satisfying
\eqref{eq_cond_mutliplier} at each $x$ with a constant $C$ independent of $x$.
\end{itemize}

The condition on the commutation with the spectral measure of $\pi(\cR)$
seems to be reasonable: in the corresponding version of the sharp
G{\aa}rding inequality on compact Lie groups in \cite{RTg} or in the abelian case, this condition
is  automatically satisfied there because if $\cR=-\cL$ and $\cL$ is the
Laplacian then $\pi(\cR)$ is central.

In \cite{HN}, Helffer and Nourigat proved that $\cR$ and, equivalently, $I+{\cR}$ are hypoelliptic. 
We finally give the result stating that the Schwartz version of such hypoellipticity
is also true.
\medskip
 
\begin{theorem}\label{thm:3}
 The operator $I+\cR$ is Schwartz-hypoelliptic, i.e. for $f\in\cS'(G)$,
 the condition $(I+\cR)f\in\cS(G)$ implies $f\in\cS(G).$
\end{theorem}
\medskip

In fact, our symbolic calculus allows the construction of a parametrix for the 
operator $I+\cR$ from which both the hypoellipticity and the Schwartz hypoellipticity follow.





\section*{Acknowledgements}
The first author acknowledges the support of the London Mathematical Society via the Grace Chisholm Fellowship. 
It was during this fellowship held at King's College London in 2011 that the work was initiated.
The second author was supported in part by the EPSRC Leadership Fellowship  EP/G007233/1.

\end{document}